\documentclass[12pt,reqno]{amsart}

\usepackage{amscd, amsfonts, amsmath, amssymb, amsthm, mathrsfs, appendix, enumerate, graphicx, flexisym, epsfig, float, tabu, graphics, hyperref, pdfpages, tikz, mathtools, comment, tikz-cd}

\usepackage{cite}

\usepackage[all,cmtip]{xy}
\usepackage{listings}

\makeatletter
\@namedef{subjclassname@2020}{Mathematics Subject Classification \textup{2020}}

\newsavebox{\@brx}
\newcommand{\llangle}[1][]{\savebox{\@brx}{\(\m@th{#1\langle}\)}%
  \mathopen{\copy\@brx\kern-0.5\wd\@brx\usebox{\@brx}}}
\newcommand{\rrangle}[1][]{\savebox{\@brx}{\(\m@th{#1\rangle}\)}%
  \mathclose{\copy\@brx\kern-0.5\wd\@brx\usebox{\@brx}}}

\theoremstyle{definition}

\numberwithin{equation}{subsection}

\usepackage[all,cmtip]{xy}

\newtheorem{thm}{Theorem}[section]

\newtheorem{cor}[thm]{Corollary}

\newtheorem{ex}[thm]{Example}

\title[Non-inner automorphisms of finite $p$-groups]{On the non-inner automorphism conjecture of finite $p$-groups}

\author{Sandeep Singh}
\address{Department of Mathematics, Akal University, Talwandi Sabo, Punjab 151302, India.}
\email{sandeepinsan86@gmail.com}

\author{Hemant Kalra}
\address{Department of Mathematics, Guru Jambheshwar University of Science and Technology, Hisar, Haryana 125001, India.}
\email{happykalra26@gmail.com}

\author{Rohit Garg}
\address{School of Mathematical Sciences, National Institute of Science Education and Research, Bhubaneswar, HBNI, P.O. Jatni, Khurda, Odisha 752050, India.}
\email{rohitgarg289@gmail.com, rohitgarg289@niser.ac.in}

\subjclass[2020]{Primary 20D15; Secondary 20D45}
\keywords{Finite $p$-group, Non-inner automorphism.}

\begin{document}
\begin{abstract}
\noindent A long-standing conjecture asserts that every finite non-abelian $p$-group
has a non-inner automorphism of order $p$.  In this paper, we settle the conjecture for a finite $p$-group ($p >2$) of nilpotency class $n$  with certain conditions. 
\end{abstract}

\maketitle

\section{Introduction}

There is a famous conjecture known as the Non-inner Automorphism Conjecture, listed in the renowned book ``Unsolved Problems in Group Theory: The Kourovka Notebook", which states that {\it Every finite non-abelian $p$-group admits an automorphism of order $p$ which
is not an inner.} (see \cite[Problem 4.13]{khu})

\medskip

Some researchers showed interest in proving the sharpened version of the conjecture. They were interested in proving that every finite non-abelian $p$-group $G$  has a non-inner automorphism of order $p$ which fixes $\Phi(G)$ element-wise (for instance, see \cite{dea, gar,gho2013,gho2014,gho2018}, and \cite{gho2021} for other references). The conjecture was first attacked by Liebeck \cite{lie}. He proved that for an odd prime $p$, every finite $p$-group  $G$ of nilpotency class $2$ has a non-inner automorphism of order $p$ fixing $\Phi(G)$ element-wise. In 2013, Abdollahi et al. \cite{abd2013b} proved the validity of the conjecture for finite $p$-groups of nilpotency class $3$. In particular, in Theorem 4.4, they proved that every finite $p$-group $G$ of odd order and of nilpotency class 3 has a non-inner automorphism of order $p$ that fixes $\Phi(G)$ element-wise. In 2014, Abdollahi et al. \cite{abd2014} showed that every finite $p$-group $G$ of co-class 2 has a non-inner automorphism of order $p$ leaving $Z(G)$ element-wise fixed. In 2017, Ruscitti et al. \cite{rus} confirmed the conjecture for finite $p$-groups of co-class 3, with $p \neq 3$.

\medskip

If there is a maximal subgroup $M$ of a finite $p$-group $G$ with $|G| > p$ and $Z(M) \subseteq Z(G)$, then there exists a non-inner automorphism of $G$ of order $p$ (see, Rotman \cite[Lemma 9.108]{rot}).  In 2002, Deaconescu and Silberberg \cite{dea} proved that if the conjecture is false for a finite $p$-group $G$, then $Z(G) < Z(M)$ for all maximal subgroups $M$ of $G$.   This raises the following natural question: 

\medskip

{\bf Question.} {\it Given a finite $p$-group $G$ with $Z(G)<Z(M)$ for all maximal subgroups $M$ of $G$, does the conjecture hold?}

\medskip

In Theorem 2.1, we  prove that every finite $p$-group $G$, ($p>2$) of nilpotency class $n$ such that $\mbox{exp}(\gamma_{n-1}(G))=p$, $|\gamma_n(G)|=p$ and $Z(C_G(x)) \le \gamma_{n-1}(G)$ for all $x \in \gamma_{n-1}(G) \setminus Z(G)$, has a non-inner automorphism of order $p$ which fixes $\Phi(G)$ element-wise. As a consequence, in Corollaries 2.2 and 2.3, we give an affirmative answer to the above question under some conditions. In \cite{rus}, the authors proved the conjecture for all non-abelian finite $p$-groups of co-class 3, where $p$ is a prime number such that $p \ne 3$. We also validate the conjecture for some non-abelian finite 3-groups of co-class 3 in Corollary 2.4.

\medskip

Throughout $p$ denotes an odd prime number. For a group $G$, by $Z_m(G)$, $\gamma_{m}(G)$, $d(G)$ and $\Phi(G)$, we denote, the $m$th term of the {\it upper central series} of $G$, the $m$th term of the {\it lower central series} of $G$, the {\it minimum number of generators} of $G$ and the {\it Frattini subgroup} of $G$, respectively. The {\it nilpotency class} and the {\it exponent} of a finite group $G$ are denoted by $\mbox{cl}(G)$ and $\mbox{exp}(G)$, respectively. A finite $p$-group $G$ of order $p^n$ with $\mbox{cl}(G)=n-c$ is said to be of {\it co-class} $c$. All other unexplained notations, if any, are standard.

\section{Main results}

Since the conjecture is true for all finite $p$-groups $G$ having nilpotency class 2 and 3, we consider only finite $p$-groups $G$ with $\mbox{cl}(G) \ge 4$.

\medskip

\begin{thm}
Let $G$ be a finite $p$-group ($p>2$) of class $n$ such that $|\gamma_{n}(G)|=\mbox{exp}(\gamma_{n-1}(G))=p$ and $Z(C_G(x)) \le \gamma_{n-1}(G)$ for all $x \in \gamma_{n-1}(G) \setminus Z(G)$. Then $G$ has a non-inner automorphism of order $p$ that fixes $\Phi(G)$ element-wise. 
\end{thm}

\begin{proof}
Since $n=\mbox{cl}(G) \ge 4$ and $\mbox{exp}(\gamma_{n-1}(G))=p$,  there exists an element $x \in \gamma_{n-1}(G)\setminus Z(G)$  of order $p$. Thus $[x,G]\subseteq \gamma_n(G),$ and therefore the order of conjugacy class of $x$ in $G$ is $p$. It follows that  $M=C_{G}(x)$ is a maximal subgroup of $G$. Let $g\in G\setminus M$. Then 
\[
(gx)^{p}=g^{p}x^{p}[x,g]^{p(p-1)/2} =g^{p}.
\]
Consider the map $\beta$ of $G$ defined as $\beta(g)=gx$ and $\beta(m)=m$ for all $m\in M$. The map $\beta$ can be extended to an automorphism of $G$ fixing $\Phi(G)$ element-wise and of order $p$. We claim that $\beta$ is a non-inner automorphism of $G$. For a contradiction, assume that $\beta=\theta_y$, the inner automorphism of $G$ induced by some $y \in G$, which implies that $y\in C_G(M)$. If $y\notin M,$ then $G=M\langle y\rangle$. It follows that $y \in Z(G)$, which is a contradiction. Therefore $y\in Z(M)$. Since $\beta=\theta_y$, we have $g^{-1}\theta_y(g)=[g,y]=x$. Now, by the given hypothesis $Z(C_G(x)) \le \gamma_{n-1}(G)$ for all $x \in \gamma_{n-1}(G) \setminus Z(G)$, we have  $y \in \gamma_{n-1}(G)$. Therefore
\[
x=[g,y] \in \gamma_n(G) \le Z(G),
\]
which contradicts  our  choice of $x$ in $G$. Hence $G$ has a non-inner automorphism of order $p$ that fixes $\Phi(G)$ element-wise. 
\end{proof}

\vspace{10 pt}

Let $G$ be a finite $p$-group such that $|Z(G)|=p$. Let $M$ be any maximal subgroup of $G$. Since $Z(M)$ is a characteristic subgroup of $M$ and $M$ is a normal subgroup of $G$, we have $Z(M)$ is a normal subgroup of $G$. Thus $Z(G)\le Z(M)$ for all maximal subgroups $M$ of  $G$. Hence, we obtain the following Corollary from Theorem 2.1:

\medskip

\begin{cor}
Let $G$ be a finite $p$-group ($p>2$) of class $n$ such that $|Z(G)|=\mbox{exp}(\gamma_{n-1}(G))=p$ and $Z(C_G(x)) \le \gamma_{n-1}(G)$ for all $x \in \gamma_{n-1}(G) \setminus Z(G)$. Then $G$ has a non-inner automorphism of order $p$ that fixes $\Phi(G)$ element-wise. 
\end{cor}

\medskip

\begin{cor}
Let $G$ be a finite $p$-group ($p>2$) of class $n$ such that $|Z(G)|=p$ and $Z(M)= \gamma_{n-1}(G)$ is of exponent $p$ for all maximal subgroups $M$ of $G$. Then $G$ has a non-inner automorphism of order $p$ that fixes $\Phi(G)$ element-wise.
\end{cor}

\begin{proof}

Given that $\mbox{cl}(G)=n$. It follows that $\gamma_n(G) \le Z(G)$. Consequently, $|\gamma_n(G)|=p$. Considering the provided hypothesis $Z(M)= \gamma_{n-1}(G)$ is of exponent $p$ for all maximal subgroups $M$ of $G$ and the proof of Theorem 2.1, we deduce that $Z(C_G(x)) \le \gamma_{n-1}(G)$ for all $x \in \gamma_{n-1}(G) \setminus Z(G)$. Hence, by Theorem 2.1, $G$ possesses a non-inner automorphism of order $p$ that fixes $\Phi(G)$ element-wise.
\end{proof}


\medskip

\begin{cor}
Let $G$ be a finite 3-group of order $3^n$ and of co-class 3 such that $Z(M)=\gamma_{n-4}(G)$ is of exponent 3 for all maximal subgroups $M$ of $G$. Then $G$ has a non-inner automorphism of order 3.
\end{cor}

\begin{proof}
Assume that $G$ does not possess any non-inner automorphism of order 3. Then, it follows from \cite[Corollary 2.3]{abd2010} that
\begin{equation}
d(Z_2(G)/Z(G))=d(Z(G)) ~ d(G).\tag{1}
\end{equation}

Since $G$ is of co-class 3, we have  $p^i \le |Z_i(G)| \le p^{i+2}$ and $\gamma_{n-i-2}(G) \le Z_i(G)$ for all $1 \le i \le n-4$. Thus, by equation (1), $d(Z(G))=1$. Now, if $|Z(G)|=p^3$, then $|Z_2(G)|=p^4$, which contradicts equation (1). Furthermore, $|Z(G)|$ cannot be $p^2$ according to \cite[Theorem 4.3]{rus}. Finally, assume that $|Z(G)|=p$. In this case, the conclusion follows from Corollary 2.3.
\end{proof}


\vspace{10 pt}

We conclude the paper by giving an example of a 3-group of order $3^7$ which supports Theorem 2.1.

\vspace{10 pt}

\begin{ex}
Consider the following group:
\[
G=\langle f_1, f_2, f_3, f_4, f_5, f_6, f_7 \rangle,
\]
with relations:
$f_3=[f_2,f_1],
f_4=f_1^3,
f_5=[f_3,f_1], 
f_6=[f_3,f_2],
f_7=[f_5,f_1],
f_7^2=[f_4,f_2], 
f_2^3=f_3^3=f_4^3=f_5^3=f_6^3=f_7^3=
[f_4,f_1]=[f_6,f_1]=[f_7,f_1]=
[f_5,f_2]=[f_6,f_2]=[f_7,f_2]=
[f_4,f_3]=[f_5,f_3]=[f_6,f_3]=[f_7,f_3]=
[f_5,f_4]=[f_6,f_4]=[f_7,f_4]=
[f_6,f_5]=[f_7,f_5]=
[f_7,f_6]=1$. Then

\medskip

\begin{enumerate}
\item[$\bullet$] $|G|=3^7$.
\vspace{5 pt}
\item[$\bullet$] The nilpotency class of $G$ is 4.
\vspace{5 pt}
\item[$\bullet$] $Z(G)=\langle f_6,f_7 \rangle$.
\vspace{5 pt}
\item[$\bullet$] $\Phi(G)=\langle f_3, f_4, f_5, f_6, f_7 \rangle$.
\vspace{5 pt}
\item[$\bullet$] $\gamma_3(G)=\langle f_5, f_6, f_7 \rangle$.
\vspace{5 pt}
\item[$\bullet$] $\gamma_4(G)=\langle f_7 \rangle$.
\end{enumerate}

\medskip

Let $x=f_5$ and $M=C_G(x)$. Then $x \in \gamma_3(G) \setminus Z(G)$ is of order 3 and $Z(M)=\langle f_5, f_6,f_7 \rangle=\gamma_3(G)$. Consider the following automorphism:
\[
\begin{array}{lcr}
\alpha(f_1f_4^2f_5f_6^2f_7^2)
=
f_1f_4^2f_5^2f_6^2f_7^2, \;
\alpha(f_1f_3f_5f_6^2)
=
f_1f_3f_5^2f_6^2, \; 
\alpha(f_1^2f_2^2f_3f_4^2f_6)
=
f_1^2f_2^2f_3f_4^2f_5^2f_6f_7.
\end{array}
\]
Now, by using the relators of $G$, we have $\alpha(f_i)=f_i ~\mbox{for all}~ i \in \{ 2,3,4,5,6,7 \} ~\mbox{and}~ \alpha(f_1)=f_1f_5$. It is easy to verify that $\alpha$ is a non-inner automorphism of order 3 which fixes $\Phi(G)$ element-wise.
\end{ex}

\section{Acknowledgement}
The research of the first author is supported by SERB, Department of Science and Technology,
under grant MTR/2022/000331. The third author expresses deep thanks to the National Institute of Science Education and Research, Bhubaneswar, India for supporting the post-doctoral research.

\section{Conflict of Interest}

The authors declare that they have no conflict of interest.

\end{document}